\documentclass[11pt,a4paper]{article}
\usepackage{lscape}
\usepackage[centertags]{amsmath}
\usepackage{amsfonts}
\usepackage{amssymb}
\usepackage{amsthm}
\usepackage{newlfont}
\usepackage{lscape}

\newtheorem{lemma}{\sc \bf Lemma}[section]
\newtheorem{propos}{\sc \bf Proposition}[section]
\newtheorem{theor}{\sc \bf Theorem}[section]
\newtheorem{corr}{\sc \bf Corollary}[section]
\newtheorem{remark}{\sc \bf Remark}[section]
\newtheorem{definition}{\sc \bf Definition}[section]
\newtheorem{example}{\sc \bf Example}[section]

\def\bkR{{\rm I\kern-.17em R}}

\begin{document}

\centerline{} \centerline{} \centerline{\Large{\bf An Example of
the Curvature Tensor for a Quantum Space}}
\centerline{}\centerline{\Large{\bf
 }}
\centerline{\bf{ Vida Milani }}
 \centerline{  Department of Mathematics, Faculty of math. Sciences, Shahid Beheshti university, Tehran, Iran}
 \centerline{e-mail: v-milani@cc.sbu.ac.ir}
  \centerline{\bf{Seyed M.H. Mansourbeigi}}
 \centerline{Department of Electrical Engineering, Polytechnic
 University, NY, USA}
 \centerline{e-mail: s.mansourbeigi@ieee.org}
 \centerline{\bf{Farzaneh Falahati  }}
\centerline{Department of Mathematics, Science and Research
Branch,Islamic Azad University, Tehran, Iran}

%
%
\begin{abstract}

\end{abstract}
{\small The paper is constructed in two parts. In the first part
we introduce the concept of the algebra $\mathbb{R}_{Q}^{2}$ of
Q-meromorphic functions on the quantum plane. The
$A_{1}(q)$-algebra of Q-analytic functions considered in [6] is
seen as a proper subalgebra.\\
In the second part we find a formula for the curvature tensor on
 $\mathbb{R}_{Q}^{2}$.\\
It is seen that when the quantization parameter tends to 1, then
this formula gives the flatness of the usual $\mathbb{R}^{2}$.
\textit{}\\

\section{Introduction}
Non-commutative geometry and quantum groups are applied to
problems of physics in different ways. Classical and quantum
mechanics on the Manin quantum plane have
been studied in [1,4,5].\\
A non-commutative framework for calculus and differential
geometry from point of view of discrete calculus has been
introduced by Kauffman [2]. In his work he modeled the notions of
derivations and derivatives with respect to different parameters
by commutators. \\
Our main objection is to transfer classical mechanics on a Poisson
algebra to its functional quantization in the sence of
[6].\\
More precisely we want to define an analogue of the poisson
bracket on $\mathbb{R}_{Q}^{2}$ and  develop an appropriate
classical mechanics on $\mathbb{R}_{Q}^{2}$ parallel to that on
$\mathbb{R}^2$.\\
On the first step, following the new interpretation of the Manin
quantum plane in [6]and by a completion of [7], we provide the
necessary tools which enable us to introduce the notions of
covariant derivative and curvature tensor on $\mathbb{R}_{Q}^{2}$.\\
In order that we can apply Kauffman modeling of derivations we
have to deal with a sufficiently large class of functions. This is
done in part two of this paper, where the formula for the
curvature tensor on $\mathbb{R}_{Q}^{2}$ is obtained showing that
when the quantization parameter goes to 1, this give us the
flateness of $\mathbb{R}^{2}$. Also a generalized formula for the
poisson bracket of two elements in $\mathbb{R}_{Q}^{2}$ is given
and properties are studied.
\section{ The algebra of  the Q-Meromorphic Functions}
Let D=$\{ q\in \mathbb{C}:|q|\leq 1\}$ be the unit unit disc in
$\mathbb{C}$. Recalling from [6], $A_{1}(q)$ will be the
$\mathbb{C}$-algebra of the all absolutely convergent power series
$\sum_{i=0}^{\infty}a_{i}q^{i}$ in D with values in $\mathbb{C}$.
Also we denote by $A_{0}(q)$ the $\mathbb{C}$-algebra of all
absolutely convergent power series $\sum_{i>-\infty}c_{i}q^{i}$ on
D-$\{0\}$ with values in $\mathbb{C}$. We can generalize the
concept of $Q$-analytic functions on the 2-intervals of
$\mathbb{R}^{2}$ with values in $A_{1}(q)$ to the algebra of
$Q$-analytic functions on the 2-interval of $\mathbb{R}^{2}$ with
values in $A_{0}(q)$ without any difficultly. Assume now that
$\Omega=\mathbb{R}-\{0\}\times \mathbb{R}-\{0\}$ and let
\begin{center}
$f=\sum_{i>-\infty}\sum_{j,k>>-\infty}a_{ijk}q^{i}t_{1}^{j}t_{2}^{k}~~~~~~~~~~~~~(2-1)$
\end{center}
be an absolutely convergent power series on $D-\{0\}\times\Omega$
with values in $\mathbb{C}$.(The sign $>>$ under the second
$\sum$ indicates j, k are bounded below). Clearly we can consider
f as a function from $\Omega$ into $A_{0}(q)$ admitting the
absolutely convergent Laurent expansion

\begin{center}
$f=\sum_{i,j>>-\infty}a_{ij}(q)t_{1}^{i}t_{2}^{j}~~~~~~~~~~~~~~~~~~~~~~~~(2-2)$
\end{center}
on $\Omega$. Since the above series is absolutely convergent on
$\Omega$, we can also write it as
\begin{center}
$f=\displaystyle\sum_{i,j=0}^{\infty\atop
-}t_{1}^{-i}\alpha_{ij}(t_{1},
t_{2})t_{2}^{-j}~~~~~~~~~~~~~~~~~~~~~~~~~(2-3)$
\end{center}
where the $\alpha_{i,j}s$ are absolutely convergent power series
on $\mathbb{R}^{2}$ with values in $A_{0}(q)$ and the sing - over
the $\sum$
means that the indices are bounded above.\\
\textbf{\textit{Definition 2-1: }} With the above notations and
conventions let
\begin{center}
$\hat{f}=\displaystyle\sum_{i,j=0}^{\infty\atop
-}x^{-i}\hat{\alpha}_{ij}(x,p)p^{-j}~~~~~~~~~~~~~~~~~~~~~~~~~~(2-4)$
\end{center}

be obtained from $f$ by the correspondence

\begin{center}
$t_{1}^{i}t_{2}^{j}=t_{2}^{j}t_{1}^{i}\longrightarrow x^{i}p^{j}$.
\end{center}
We call $\hat{f}$ a Q-meromorphic function on $\Omega$ with
values in  $A_{0}(q)$ or simply a Q-meromorphic  function on
$\Omega$.\\
The two functions $\frac{1}{x}$ and  $\frac{1}{p}$ are
Q-meromorphic functions on $\Omega$ satisfying the following
commutation relations
\begin{center}
$x\frac{1}{x}=\frac{1}{x}x=1,~ p\frac{1}{p}=\frac{1}{p}p=1~~~~~~~~~~~~~~~~~~~~~~~~~~~~(2-5)$\\
$p\frac{1}{x}=q^{-1}\frac{1}{x}p,~ x\frac{1}{p}=q\frac{1}{p}x~~~~~~~~~~~~~~~~~~~~~~~~~~~~~~~~~(2-6)$\\
$\frac{1}{p}\frac{1}{x}=q\frac{1}{x}.\frac{1}{p},~(\frac{1}{x})^{i}=\frac{1}{x^{i}},~(\frac{1}{p})^{j}=\frac{1}{p^{j}}~~~~~~~~~~~~~~~~~~~~~~(2-7)$.
\end{center}

By using these commutation relations we always follow the order
$(x^{i}p^{j})_{i,j>>-\infty}$
 in writing the Q-meromorphic functions as above.\\

\textbf{\textit{Remark 2-1: }}If $\hat{f}(x,p)$ is a Q-analytic on
 $\Omega$ with values in $A_{0}(q)$, then for k,l$\in$ $\mathbb{Z}$,
$\hat{f}(q^{k}x,q^{l}p)$ is also a Q-analytic function on $\Omega$ with values in $A_{0}(q)$. \\

\textbf{\textit{Definition 2-2: }} The product of two
Q-meromorphic functions
\begin{center}
$\hat{f}=\displaystyle\sum_{i_{1},i_{2}=0}^{\infty\atop -}x^{-i_{1}}\hat{a}_{i_{1}i_{2}}(x,p)p^{-i_{2}}$\\
$\hat{g}=\displaystyle\sum_{j_{1},j_{2}=0}^{\infty\atop
-}x^{-j_{1}}\hat{b}_{j_{1}j_{2}}(x,p)p^{-j_{2}}$
\end{center}
on $\Omega$ will be defined by \\

$\hat{f}.\hat{g}=\displaystyle\sum_{i_{1},i_{2}=0}^{\infty\atop
-}\sum_{j_{1},j_{2}=0}^{\infty\atop
-}q^{i_{2}j_{1}}x^{-i_{1}-j_{1}}(\hat{a}_{i_{1}i_{2}}(x,q^{-j_{1}}p).\hat{b}_{j_{1}j_{2}}(q^{-i_{2}}x,p))p^{-i_{2}-j_{2}}~~~~(2-8)$

where the above product between $\hat{a}_{i_{1}i_{2}}$ and
$\hat{b}_{j_{1}j_{2}}$ is the product of two Q-analytic functions
with values in $A_{0}(q)$ in the sense of
[6].\\

\textbf{\textit{Lemma 2-1: }} With the above notations the product
of two Q-meromorphic functions $\hat{f}$ and $\hat{g}$ on is
Q-meromophic function on $\Omega$\\

 \textbf{\textit{proof. }}The proof is easily seen from the fact
 that
 $\hat{a}_{i_{1}i_{2}}(x,q^{-j_{1}}p).\hat{b}_{j_{1}j_{2}}(q^{-i_{2}}x,p)$
 is a Q-analytic function on the quantum plane with values in
 $A_{0}(q)$.\\

From the above lemma we can see that the set of all Q-meromorphic
functions on $\Omega$ with values in $A_{0}(q)$ is a
non-commutative, associative, unital $A_{0}(q)$-algebra. This
algebra  which we denote here after  by $\mathbb{R}_{Q}^{2}$,
contains $A_{Q}$, the $A_{1}(q)$-algebra of Q-analytic functions
on the quantum plane with values in A$_{1}$(q), as its subalgebra.
It is clear that $\mathbb{R}_{Q}^{2}$ is the (1, D-0, $A_{0}(q)$)
functional quantization of $M$: the $\mathbb{C}$-algebra of all
absolutely convergent power series
$\sum_{i,j>>-\infty}a_{ij}t_{1}^{i}t_{2}^{j}$ on $\Omega$ with
values in $\mathbb{C}$ in the sense of [6], and if we denote by
$A$ the $\mathbb{C}$-algebra of all entire functions of the form
$\sum_{i,j=0}^{\infty}a_{ij}t_{1}^{i}t_{2}^{j}$ on
$\mathbb{R}^{2}$ with
values in $\mathbb{C}$, then $i_A o \Phi_{A} = \Phi_{M} o i_{A_Q}$\\

Where $\Phi_{A}$ and $\Phi_{M}$ are the quantization maps defined
in[6], $A_{Q}\longrightarrow \mathbb{R}_{Q}^{2}$ and
$A\longrightarrow M$ are the canonical injections and $i_A : A
\rightarrow M$ and $i_{A_Q} : A_Q \rightarrow \mathbb{R}_{Q}^{2}$
are the inclusions.

\section{ Derivative and the curvature tensor }
Following Wess and Zumino [3,8]  we can generalize the
differential calculus by defining differential operators $\partial
x$ and
$\partial y$ as follows:\\
 $\frac{\partial}{\partial
x}x=1+q^{2}x\frac{\partial}{\partial
x}+(q^{i}-1)y\frac{\partial}{\partial
y}~~~~~~~~~~~~~~~~~~~~~~~~~~~~~~~~~~~~~~~~~~~~~~~~~~~~~~~~~~~~~~~~~~~~~~~~(3-1)$\\\\
$\frac{\partial}{\partial x}y=qy\frac{\partial}{\partial
x}~~~~~~~~~~~~~~~~~~~~~~~~~~~~~~~~~~~~~~~~~~~~~~~~~~~~~~~~~~~~~~~~~~~~~~~~~~~~~~~~~~~~~~~~~~~~~~~~~~~~~~~~~~~~~~~~~~~~~~~(3-2)$\\\\
$\frac{\partial}{\partial y}x=q^{-1}x\frac{\partial}{\partial
y}~~~~~~~~~~~~~~~~~~~~~~~~~~~~~~~~~~~~~~~~~~~~~~~~~~~~~~~~~~~~~~~~~~~~~~~~~~~~~~~~~~~~~~~~~~~~~~~~~~~~~~~~~~~~~~~~~~~~~~~(3-3)$\\\\
$\frac{\partial}{\partial y}y=1+q^{2}y\frac{\partial}{\partial
y}~~~~~~~~~~~~~~~~~~~~~~~~~~~~~~~~~~~~~~~~~~~~~~~~~~~~~~~~~~~~~~~~~~~~~~~~~~~~~~~~~~~~~~~~~~~~~~~~~~~~~~~(3-4)$\\\\
from above relations, it's easy to see that: \\

$\frac{\partial}{\partial
y}(y^{n}x^{m})=y^{n-1}x^{m}\frac{1-q^{2n}}{1-q^{2}}~~~~~~~~~~~~~~~~~~~~~~~~~~~~~~~~~~~~~~~~~~~~~~~~~~~~~~~~~~~~~~~~~~~~~~~~~~~~~~(3-5)$\\

\textbf{\textit{ Definition 3-1:}}
 For a fixed element $H\in
 \mathbb{R}_{Q}^{2}$, the derivative of an arbitrary element $f\in
 \mathbb{R}_{Q}^{2}$ with respect to the time parameter t is
 defined by $\frac{df}{dt}:=[f,H]$.
From this definition it is seen that H is independent of t. \\
(for a mechanical system, H can be considered as the Hamiltonian
function of the system.)\\

In classical differential geometry the Levi-civita connection on
$\mathbb{R}^{2}$ gives us the following covariant derivatives on
functions on $\mathbb{R}^{2}$:\\

$\nabla_{x}f:=\frac{\partial f}{\partial
x}=\partial_{x}f~~~~~~~~~~~~~~~~~~~~~~~~~~~~~~~~~~~~~~~~~~~~~~~~~~~~~~~~~~~~~~~~~~~~~~~~~~~~~~~~~~~~~~~~~~~~~~~~~~~~~~~~~~~~~~~~~~~~(3-6)$\\\\

$\nabla_{y}f:=\frac{\partial f}{\partial
y}=\partial_{y}f~~~~~~~~~~~~~~~~~~~~~~~~~~~~~~~~~~~~~~~~~~~~~~~~~~~~~~~~~~~~~~~~~~~~~~~~~~~~~~~~~~~~~~~~~~~~~~~~~~~~~~~~~~(3-7)$\\\\

 \textbf{\textit{ Definition 3-2:}} for each $f\in
\mathbb{R}_{Q}^{2}$ we define the following covariant
derivatives with respect to $x$ and $y$ by:\\
$\nabla_{x}f:=[f,
H].\frac{1}{[x,H]}~~~~~~~~~~~~~~~~~~~~~~~~~~~~~~~~~~~~~~~~~~~~~~~~~~~~~~~~~~~~~~~~~~~~~~~~~~~~~~~~~~~~~~~~~~~~~~~~~~~~~~~~~~~~~~~~~(3-8)$\\\\
$\nabla_{y}f:=[f,
H].\frac{1}{[y,H]}~~~~~~~~~~~~~~~~~~~~~~~~~~~~~~~~~~~~~~~~~~~~~~~~~~~~~~~~~~~~~~~~~~~~~~~~~~~~~~~~~~~~~~~~~~~~~~~~~~~~~~~~~~~~~~~~~~~~~~~~~~(3-9)$\\\\
$\nabla_{xy}f:=[f,
H].\frac{1}{[xy,H]}~~~~~~~~~~~~~~~~~~~~~~~~~~~~~~~~~~~~~~~~~~~~~~~~~~~~~~~~~~~~~~~~~~~~~~~~~~~~~~~~~~~~~~~~~~~~~~~~~~~~~~~~~~~~~~~~~~~~~~~~~(3-10)$\\\\

\textbf{\textit{ Proposition 3-1:}} the following properties of
the covariant derivative is obtained:\\

a) $\nabla_{\lambda x}f=\frac{1}{\lambda}\nabla{x}f$\\

b) $\nabla_{[x,y]}f=(1-q)^{-1}\nabla_{xy}f
$\\

c)
$\frac{1}{\nabla_{xy}f}=\frac{x}{\nabla_{y}f}+\frac{y}{\nabla_{x}f}$\\

\textbf{\textit{proof}}:\\

a) It's obvious by using definition 3-2.\\

b) Now, we compute $\nabla_{[x,y]}f$ as follow:\\

$\nabla_{[x,y]}f=\nabla_{xy-yx}f=\nabla_{xy-qyx}f=\nabla_{(1-q)xy}f$\\\\
$=(1-q)^{-1}\nabla_{xy}f$\\

c)$\nabla_{xy}f=[f,H].\frac{1}{x[y,H]+[x,H]y}$\\\\
$\Longrightarrow
\frac{1}{\nabla_{xy}f}=(x[y,H]+[x,H]y)\frac{1}{[f,H]}$\\\\
$x\frac{1}{\nabla_{y}f}+\frac{1}{\nabla_{x}f}*y$\\

And by abuse of notation we can write it as

$\frac{1}{\nabla_{xy}f}=\frac{x}{\nabla_{y}f}+\frac{y}{\nabla_{x}f}$\\

\textbf{\textit{ Definition 3-3:}}The curvature  tensor of
$\mathbb{R}_{Q}^{2}$ is defined
as follows:\\

$R(x,y) : R^{2}_{Q}\longrightarrow \Bbb R$\\

$R(x, y)f=[\nabla _{x},\nabla_{
y}]f-\nabla_{[x,y]}f~~~~~~~~~~~~~~~~~~~~~~~~~~~~~~~~~~~~~~~~~~~~~~~~~~~~~~~~~~~~~~~~~~~~~~~~~~~~~~~~~~~~~~~~~~~~~~~~~~~~~~(3-11)$\\
Now, we compute both $[\nabla_{ x},\nabla _{y}]f$ and
$\nabla_{[x,y]}$
separately: \\
from (3-6) and (3-7), it follows that:\\
$[\nabla_{ x},\nabla_{ y}]=[\partial_{ x} , \partial_{
y}]=\partial_{ x}\partial_{y}
-\partial_{ y}\partial_{ x}$\\\\
$=\partial_{x}\partial_{y}-q\partial_{x}\partial_{y}
$\\\\
$=(1-q)\partial_{x}\partial_{y}
~~~~~~~~~~~~~~~~~~~~~~~~~~~~~~~~~~~~~~~~~~~~~~~~~~~~~~~~~~~~~~~~~~~~~~~~~~~~~~~~~~~~~~~~~~~~~~~~~~~~~~~(3-12)$\\
Now, let $f\in \mathbb{R}_{Q}^{2}$ and $f =\displaystyle\sum
_{i,j} a_{ij}(q)x^{i}y^{j}$ then\\

$[\partial_{ x}, \partial_{ y}]f=(1-q)\partial _{x}\partial_{
y}(\displaystyle\sum _{i,j} a_{ij}(q)x^{i}y^{j})$\\

$(1-q)\displaystyle\sum _{i,j} a_{ij}(q)\partial _{x}\partial_{ y}
(x^{i}y^{j})~~~~~~~~~~~~~~~~~~~~~~~~~~~~~~~~~~~~~~~~~~~~~~~~~~~~~~~~~~~~~~~~~~~~~~~~~~~~~~~~~~~~~~~~~~~~~~~~~(3-13)$ \\

$[\partial _{x},
\partial_{y}]f=\frac{1}{(1+q)^{2}(1-q)}\displaystyle\sum
_{i,j}q^{ij+j-1}(1-q^{2j})(1-q^{2i})a_{ij}(q)y^{j-1}x^{i-1}~~~~~~~~~~~~~~~~~~~~~~~~~~~(3-14)$\\

In this way we obtain the following formula for the curvature of
$\mathbb{R}_{Q}^{2}$

 $R(x,y)f=\frac{1}{(1+q)^{2}(1-q)}\displaystyle\sum
_{i,j}q^{ij+j-1}(1-q^{2j})(1-q^{2i})a_{ij}(q)y^{j-1}x^{i-1}-(1-q)^{-1}(\frac{1}{\frac{x}{\nabla_{y}f}+\frac{y}{\nabla_{x}f}})$\\

 \end{document}